\documentclass[a4paper, 12pt]{article}
\usepackage{amssymb}
\usepackage{amstext}
\usepackage{amsmath}
\usepackage{amscd}
\usepackage{amsthm}
\usepackage{latexsym}
\usepackage{amsfonts}
\usepackage[all]{xy}
\usepackage{enumerate}
\usepackage{textcomp}
\usepackage{color}
\usepackage{mathrsfs}
\usepackage{bbm}

\numberwithin{equation}{section}

\theoremstyle{plain}
\newtheorem{Theorem}{Theorem}[section]

\theoremstyle{definition}

\newtheorem{Definition}[Theorem]{Definition}

\def\HH{{\bf H}}
\def\NN{{\bf N}}
\def\hh{{\bf h}}
\def\t{{\bf t}}
\def\s{{\bf s}}
\def\aa{{\bf a}}
\def\bb{{\bf b}}
\def\ii{{\bf i}}
\def\xx{{\bf x}}

\def\uu{{\bf u}}
\def\vv{{\bf v}}
\def\kk{{\bf k}}

\def\rr{{\bf r}}

\begin{document}
\pagestyle{plain}
\setlength{\baselineskip}{20pt}

\title{Asymptotic behavior for quadratic variations of non-Gaussian multiparameter Hermite random fields}
\author{T. T. Diu Tran\footnote{Universit\'{e} du Luxembourg, UR en math\'{e}matiques, Maison du Nombre, 6 avenue de la Fonte, L-4364 Esch-sur-Alzette, Luxembourg. E-mail: diu.tran@uni.lu}}
\maketitle

\begin{abstract} 
Let $(Z^{q, \HH}_\t)_{\t \in [0, 1]^d}$ denote a $d$-parameter Hermite random field of order $q \geq 1$ and self-similarity parameter $\HH = (H_1, \ldots, H_d) \in (\frac{1}{2}, 1)^d$. This process is $\HH$-self-similar, has stationary increments and exhibits long-range dependence. Particular examples include fractional Brownian motion ($q=1$, $d=1$), fractional Brownian sheet $(q=1, d \geq 2)$, Rosenblatt process ($q=2$, $d=1$) as well as Rosenblatt sheet $(q=2, d \geq 2)$. For any $q \geq 2, d\geq 1$ and $\HH \in (\frac{1}{2}, 1)^d$  we show in this paper that a proper normalization of the quadratic variation of $Z^{q, \HH}$ converges in $L^2(\Omega)$ to a standard $d$-parameter Rosenblatt random variable with self-similarity index $\HH'' = 1+ (2\HH-2)/q$.
\end{abstract}

 

\section{Motivation and main results}

In recent years, analysing the asymptotic behaviour of power variations of self-similar stochastic processes has attracted a lot of attention. This is because they play an important role in various aspects, both in probability and statistics. As far as quadratic variations are concerned, a classical application  is to use them for the construction of efficient estimators for the self-similarity parameter (see e.g. \cite{TudorViens1, TudorViens2}). For a less conventional application, let us also mention the recent reference \cite{IvanG}, in which the authors have used weighted power variations of fractional Brownian motion to compute exact rates of convergence of some approximating schemes associated to one-dimensional fractional stochastic differential equations.

In this paper, we deal with the quadratic variation in the context of  {\it multiparameter Hermite random fields}. To be more specific, let $Z^{q, \HH} = (Z^{q, \HH}_\t)_{\t \in [0, 1]^d}$ stand for the $d$-parameter Hermite random field of order $q \geq 1 $ and self-similarity parameter $\HH = (H_1,\ldots, H_d) \in (\frac{1}{2}, 1)^d$ (see Definition \ref{Deforiginal} for the precise meaning), and consider a renormalized version of its quadratic variation, namely
\begin{equation}\label{eq:1}
V_{\NN}: = \frac{1}{\NN}\sum_{\ii=0}^{\NN-1}\bigg[\NN^{2\HH}\Big(\Delta Z^{q, \HH}_{[\frac{\ii}{\NN}, \frac{\ii +1}{\NN}]}\Big)^2 - 1 \bigg],
\end{equation} 
where $\Delta Z^{q, \HH}_{[\s, \t]}$ is the increments of $Z^{q, \HH}$ defined as
\begin{equation}\label{eq:incrementZ}
\Delta Z^{q, \HH}_{[\s, \t]} = \sum_{\rr \in \{0, 1\}^d} (-1)^{d- \sum_i r_i} Z^{q, \HH}_{\s + \rr.(\t - \s)},
\end{equation}
and where the bold notation is systematically used in presence of multi-indices (we refer to Section 2 for precise definitions). As illustrating examples, observe that (\ref{eq:incrementZ}) reduces to $\Delta Z^{q, H}_{[s, t]} = Z^{q, H}_t -Z^{q, H}_s$ when $d=1$, and to  $\Delta Z^{q, H_1, H_2}_{[\s, \t]} = Z^{q, H_1, H_2}_{t_1, t_2} - Z^{q, H_1, H_2}_{t_1, s_2} - Z^{q, H_1, H_2}_{s_1, t_2} + Z^{q, H_1, H_2}_{s_1, s_2}$ when $d=2$. 

It is well-known that each Hermite random field $Z^{q, \HH}$ is $\HH$-self-similar (that is, $ (Z^{q, \HH}_{\aa \t})_{\t \in \mathbb{R}^d} \overset{(d)}{=} (\aa^\HH Z^{q, \HH}_\t)_{\t \in \mathbb{R}^d}$ for any $\aa > 0$), has stationary increments (that is $(\Delta Z^{q, \HH}_{[0, \t]})_{\t \in \mathbb{R}^d} \overset{(d)}{=} (\Delta Z^{q, \HH}_{[\hh, \hh+ \t]})_{\t \in \mathbb{R}^d}$ for all $\hh \in\mathbb{R}^d$) and exhibits long-range dependence. Also, when $q=1$, observe that $Z^{1, \HH}$ is either the fractional Brownian motion (if $d=1$) or the fractional Brownian sheet (if $d\geq 2$); in particular, among all the Hermite random fields $Z^{q,\HH}$, it is the only one to be Gaussian. When $q=2$, we use the usual terminologies Rosenblatt process (if $d=1$) or Rosenblatt sheet (if $d \geq 2$). 

Before describing our results, let us give a brief overview of the current state of the art. Firstly, let us consider the case  $q=d=1$, that is, the case where $Z^{1, H} = B^H$ is a fractional Brownian motion with Hurst parameter $H$. The behavior of the quadratic variation of $B^H$ is well-known since the eighties, and dates back to the seminal works of Breuer and Major \cite{BreuerMajor}, Dobrushin and Major \cite{DobrushinMajor}, Giraitis and Surgailis \cite{Giraitis} or Taqqu \cite{Taqqu}. We have, as $N\to\infty$:
\begin{itemize}
\item If $H<3/4$, then 
$$N^{-1/2}\sum_{j=1}^N \bigg(N^{2H}\Big(B^H_{j/N} - B^H_{(j-1)/N}\Big)^2 - 1\bigg) \xrightarrow[]{(d)} \mathcal{N}(0, \sigma_H^2).$$
\item If $H =3/4$, then
$$(N\log{N})^{-1/2}\sum_{j=1}^N \bigg(N^{3/2}\Big(B^H_{j/N} - B^H_{(j-1)/N}\Big)^2 - 1\bigg) \xrightarrow[]{(d)} \mathcal{N}(0, \sigma_{3/4}^2).$$
\item If $H >3/4$, then
$$N^{1-2H}\sum_{j=1}^N \bigg(N^{2H}\Big(B^H_{j/N} - B^H_{(j-1)/N}\Big)^2 - 1\bigg) \xrightarrow[]{L^2(\Omega)} \text{``Rosenblatt r.v''},$$
where $\text{``Rosenblatt r.v''}$ denotes the random variable which is the value at time 1 of the Rosenblatt process.
\end{itemize}
Secondly, assume now that $q=1$ and $d =2$, that is, consider the case where $Z^{1, \HH}$ is this time a two-parameter fractional Brownian sheet with Hurst parameter $\HH = (H_1, H_2)$. According to R\'eveillac, Stauch and Tudor \cite{ReveillacTudor} and with $\varphi(N, \HH)$ a suitable scaling factor, the quadratic variation of $Z^{1, \HH}$ behaves as follows, as $N\to\infty$:
\begin{itemize}
\item If $\HH\notin (3/4, 1)^2$, then
$$\varphi (N, \HH)\sum_{i=1}^N\sum_{j=1}^N \bigg(N^{2H_1 + 2H_2}\Big(\Delta Z^{1, \HH}_{[\frac{\ii-1}{\NN}, \frac{\ii }{\NN}]}\Big)^2 - 1\bigg) \xrightarrow[]{(d)} \mathcal{N}(0, \sigma_{\HH}^2).$$
\item If $\HH \in (3/4, 1)^2$, then
\begin{align*}
\varphi (N, \HH)\sum_{i=1}^N\sum_{j=1}^N \bigg(N^{2H_1 + 2H_2}&\Big(\Delta Z^{1, \HH}_{[\frac{\ii-1}{\NN}, \frac{\ii }{\NN}]}\Big)^2 - 1\bigg)\\ 
&  \xrightarrow[]{L^2(\Omega)} \text{``two-parameter Rosenblatt r.v"},
\end{align*}
where $\text{``two-parameter Rosenblatt r.v''}$ means the value at point ${\bf 1}=(1, 1)$ of the two-parameter Rosenblatt sheet. 
\end{itemize}
Here, we observe the following interesting phenomenon: the limit law in the mixture case (that is, when $H_1 \leq 3/4$ and $H_2 > 3/4$) is Gaussian. For the simplicity of exposition, above we have only described what happens when $d=2$. But the asymptotic behaviour for the quadratic variation of $Z^{1, \HH}$ is actually known for any value of the dimension $d\geq 2$, and we refer to Pakkanen and R\'{e}veillac \cite{Reveillac2, Reveillac3, Reveillac1} for precise statements.

Let us finally review the existing literature about the quadratic variation of $Z^{q,\mathbf{H}}$ in the {\it non}-Gaussian case, that is, when $q\geq 2$. It is certainly because it is a more difficult case to deal with that only the case where $d=1$ has been studied so far. Chronopoulou, Tudor and Viens have shown in \cite{TudorViens1} (see also \cite{TudorViens2, Tudor}) that, properly renormalized, the quadratic variation of $Z^{q, H}$ converges in $L^2(\Omega)$, for \textit{any} $q \geq 2$ and \textit{any} value of $H \in (1/2, 1)$, to the Rosenblatt random variable. A consequence of this finding is that fractional Brownian motion is the only Hermite process ($d=1$) for which there exists a range of parameters such that its quadratic variation exhibits {\it normal} convergence; indeed, for all the other Hermite processes, \cite{TudorViens1} shows that we have the convergence towards a non-Gaussian random variable belonging to the second Wiener chaos. 

In the present paper, we study what happens in the remaining cases, that is, when $q$ and $d$ are both bigger or equal than 2.
Thanks to our main result, Theorem \ref{main}, we now have a complete picture for the asymptotic behaviour of the quadratic variation of {\it any} Hermite random field.

\begin{Theorem} \label{main}
Fix $q \geq 2$, $d\geq 1$ and $\HH \in (\frac{1}{2}, 1)^d$. Let $Z^{q, \HH}$ be a $d$-parameter Hermite random field of order $q$ and self-similarity parameter $\HH$ (see Definition \ref{Deforiginal}). Then $c_{1,\HH}^{-1/2}\NN^{(2-2\HH)/q}(q!q)^{-1} V_{\NN}$ converges, in $L^2(\Omega)$, to the standard $d$-parameter Rosenblatt sheet with self-similarity parameter $1+ (2\HH-2)/q$ evaluated at time $\bf{1}$, where $c_{1, \HH}$ given by (\ref{eq:3.5}).
\end{Theorem}

 Our proof of Theorem \ref{main} follows a strategy introduced by Tudor and Viens in \cite{TudorViens2}, based on the use of chaotic expansion into multiple Wiener-It\^{o} integrals. Let us sketch it. Since the Hermite random field $Z^{q, \HH}$ is an element of the $q$-th Wiener chaos, we can firstly rely on the product formula for multiple integrals to obtain that the quadratic variation $V_{\NN}$ can be decomposed into a sum of multiple integrals of even orders from $2$ to $2q$, see Section 3.1. We are thus left to analyse the behavior of each chaos component. As we will prove in Section 3.2,  the dominant term of $V_\NN$ (after proper normalization) is the term in the second Wiener chaos, that is, all the other terms in the chaotic expansion are asymptotically negligible. Finally, by using the isometric property of multiple Wiener-It\^{o} integrals and checking the $L^2(([0, 1]^d)^2)$ convergence of its kernel, we will prove in Section 3.3 that the projection onto the second Wiener chaos converges in $L^2(\Omega)$ to the $d$-parameter Rosenblatt random variable, which will lead to the convergence of the normalization of $V_\NN$ to the same random variable.

In conclusion, it is worth pointing out that, irrespective of the self-similarity parameter, the (properly normalized) quadratic variation of any \textit{non-Gaussian} multiparameter Hermite random fields exhibits a convergence to a random variable belonging to the second Wiener chaos. It is in strong contrast with what happens in the Gaussian case ($q=1$), where either central or non-central limit theorems may arise, depending on the value of the self-similarity parameter.

The remainder of the paper is structured as follows. Section $2$ contains some preliminaries and useful notation. The proof of our main result, namely Theorem \ref{main}, is then provided in Section $3$.

\section{Preliminaries}

This section describes the notation and the mathematical objects (together with their main properties) that are used throughout this paper.

\subsection{Notation}

Fix an integer $d\geq 1$. In what follows,  we shall systematically use bold notation when dealing with multi-indexed quantities. We thus write $\aa = (a_1, a_2, \ldots, a_d)$, $\aa\bb = (a_1b_1, a_2b_2, \ldots, a_db_d)$ or $\aa/\bb = (a_1/b_1, a_2/b_2, \ldots, a_d/b_d)$. Similarly, $[\aa, \bb] = \prod_{i=1}^d[a_i, b_i],  (\aa, \bb) = \prod_{i=1}^d(a_i, b_i)$. Summation is as follows: $\sum_{\ii=1}^\NN a_\ii = \sum_{i_1 =1}^{N_1}\sum_{i_2 =1}^{N_2}\ldots \sum_{i_d =1}^{N_d}a_{i_1,i_2,\ldots,i_d}$ whereas, for products, we shall write $\aa^\bb = \prod_{i=1}^d a_i^{b_i}$. Finally, we shall write $ \aa < \bb$ (resp. $ \aa \leq \bb$) whenever $a_1 < b_1$, $a_2 < b_2$, $\ldots$, $a_d < b_d$
(resp. $a_1 \leq b_1$, $a_2 \leq b_2$, $\ldots$, $a_d \leq b_d$).

\subsection{Multiple Wiener-It\^{o} integrals}
We will now briefly review the theory of multiple Wiener-It\^{o} integrals with respect to the Brownian sheet, as described e.g. in Nualart's book \cite{Nualart} (chapter 1 therein) or in \cite[Section 3]{Reveillac2}. Let $f \in L^2((\mathbb{R}^{d})^q)$ and let us denote by $I_q^W(f)$ the $q$-fold multiple Wiener-It\^{o} integral of $f$ with respect to the standard two-sided Brownian sheet $(W_\t)_{\t \in \mathbb{R}^d}$. In symbols, such an integral is written
\begin{equation}\label{eq:abc}
I_q^W(f)  = \int_{(\mathbb{R}^d)^q}  dW_{\uu_1}\ldots dW_{\uu_q} f(\uu_1, \ldots, \uu_q).
\end{equation}
Moreover, one has $I_q^W(f) = I_q^W(\widetilde{f})$, where $\widetilde{f}$ is the symmetrization of $f$ defined by
\begin{equation}\label{eq:dola}
\widetilde{f}(\uu_1, \ldots, \uu_q) = \frac{1}{q!}\sum_{\sigma \in \mathfrak{S}_q} f(\uu_{\sigma(1)}, \ldots, \uu_{\sigma(q)}).
\end{equation}
The set of random variables of the form $I_q^W(f)$, when $f$ runs over $L^2((\mathbb{R}^{d})^q)$, is called the $q$th Wiener chaos of $W$. Furthermore, if $f \in L^2((\mathbb{R}^{d})^p)$ and $g \in L^2((\mathbb{R}^{d})^q)$ are two symmetric functions, then
\begin{equation}\label{eq:P1}
I_p^W(f)I_q^W(g) = \sum_{r=0}^{p \wedge  q} r!\binom{p}{r}\binom{q}{r}I_{p+q-2r}^W(f \widetilde{\otimes}_r g),
\end{equation}
where the contraction $f \otimes_r g$, which belongs to $L^2((\mathbb{R}^{d})^{p+q-2r})$ for every $r = 0, 1, \ldots, p \wedge q$, is given by
\begin{align}\label{eq:P2}
f \otimes_r g& (\uu_1, \ldots, \uu_{p-r}, \vv_1, \ldots, \vv_{q-r}) \nonumber\\ 
& = \int_{(\mathbb{R}^{d})^r} d\aa_1 \ldots d\aa_r f(\uu_1, \ldots, \uu_{p-r}, \aa_1, \ldots, \aa_r) g(\vv_1, \ldots, \vv_{q-r}, \aa_1, \ldots, \aa_r) 
\end{align}
and $f \widetilde{\otimes}_r g$ stands for the symmetrization of $f \otimes_r g$ (according to the notation introduced in (\ref{eq:dola})). 
For any $r= 0, \ldots, p\wedge q$, Cauchy-Schwarz inequality yields
\begin{equation}\label{eq:3}
\| f \widetilde{\otimes}_r g\|_{L^2((\mathbb{R}^{d})^{p+q-2r})} \leq \| f \otimes_r g\|_{L^2((\mathbb{R}^{d})^{p+q-2r})} \leq \|f\|_{L^2((\mathbb{R}^{d})^p)}\|g\|_{L^2((\mathbb{R}^{d})^{q})}.
\end{equation}
Also, $f \otimes_p  g = \left\langle f, g \right\rangle_{L^2((\mathbb{R}^{d})^p)}$ when $q=p$.
Furthermore, multiple Wiener-It\^{o} integrals satisfy the following isometry and orthogonality properties
$$E[I_p^W(f)I_q^W(g)] =
\begin{cases}
 & p! \big\langle \widetilde{f}, \widetilde{g} \big\rangle_{L^2((\mathbb{R}^{d})^{p})} \qquad\text{if } p=q\\
& 0 \qquad\qquad\qquad\quad\quad \text{if } p \ne q.
\end{cases}$$

\subsection{Multiparameter Hermite Random Fields}

Let us now introduce our main object of interest in this paper, the so-called multiparameter Hermite random field. We follow the definition given by Tudor in \cite[Chapter 4]{Tudor}.

\begin{Definition}\label{Deforiginal}
{\it
Let $q,d \geq 1$ be two integers and let $\HH = (H_1, \ldots, H_d)$ be a vector belonging to $(\frac{1}{2}, 1)^d$.
The $d$-parameter Hermite random field of order $q$ and self-similarity parameter $\HH$ means any random field of the form
\begin{align}\label{Def1}
Z^{q, \HH}(\t)&= c_{q, \HH} \int_{(\mathbb{R}^{d})^{ q}} dW_{u_{1,1}, \ldots, u_{1,d}} \ldots dW_{u_{q,1}, \ldots, u_{q,d}} \nonumber\\ 
&\qquad\times \bigg(\int_0^{t_1} da_1 \ldots \int_0^{t_d} da_d \prod_{j=1}^q (a_1-u_{j,1})_+^{-(\frac{1}{2} + \frac{1-H_1}{q})}\ldots (a_d-u_{j,d})_+^{-(\frac{1}{2} + \frac{1-H_d}{q})}\bigg)\nonumber\\
& = c_{q, \HH}\int_{(\mathbb{R}^{d})^{q}}dW_{\uu_1} \ldots dW_{\uu_q} \int_0^\t d\aa \prod_{j=1}^q (\aa - \uu_j)_+^{-(\frac{1}{2} + \frac{1-\HH}{q})},
\end{align}
where $x_+ = \max(x, 0)$, $W$ is a standard two-sided Brownian sheet, and $c(q, \HH)$ is the unique positive constant depending only on $q$ and $\HH$
chosen so that $E[Z^{q, \HH}(\textbf{1})^2] = 1$. }
\end{Definition}
The above integral  (\ref{Def1}) represents a multiple Wiener-It\^{o} integral of the form (\ref{eq:abc}).

In many occasions (for instance when one wants to simulate $Z^{q,\HH}$, or when one looks for constructing a stochastic calculus with respect to it), the following finite-time representation for $Z^{q,\HH}$ may also be of interest:
\begin{align}\label{Def2}
&Z^{q, \HH}(\t) \overset{(d)}{=} b_{q, \HH} \int_0^{t_1} \ldots \int_0^{t_d} dW_{u_{1,1}, \ldots, u_{1,d}} \ldots \int_0^{t_1}\ldots\int_0^{t_d} dW_{u_{q,1}, \ldots, u_{q,d}} \nonumber\\ 
&\hspace{2.5cm}\times \bigg(\int_{u_{1,1} \vee \ldots \vee u_{q,1}}^{t_1} da_1 \partial_1K^{H_1'}(a_1, u_{1,1}) \ldots \partial_1K^{H_1'}(a_1, u_{q,1})\bigg) \nonumber \\
& \hspace{6cm} \vdots \nonumber\\
& \hspace{2.5cm}\times \bigg(\int_{u_{1,d} \vee \ldots \vee u_{q,d}}^{t_d} da_d \partial_1K^{H_d'}(a_d, u_{1,d}) \ldots \partial_1K^{H_d'}(a_d, u_{q,d})\bigg) \nonumber\\
& = b_{q, \HH}\int_{[0, \t]^{q}}dW_{\uu_1}\ldots dW_{\uu_q} \prod_{j=1}^d\int_{u_{1,j} \vee \ldots \vee u_{q,j}}^{t_j} da \partial_1K^{H'_j}(a, u_{1,j}) \ldots \partial_1K^{H'_j}(a, u_{q,j}).
\end{align}
In (\ref{Def2}), $K^{H}$ stands for the usual kernel appearing in the
classical expression of  the fractional Brownian motion $B^H$ as a Volterra integral with respect to Brownian motion (see e.g. \cite{Ivan, Ivan1}), that is,
$B^H_t = \int_0^t K^H(t, s)dB_s,$
whereas
\begin{equation}\label{eq:bqH}
b_{q, \HH} :=\frac{(\HH(2\HH-1))^{1/2}}{(q!(\HH'(2\HH' -1))^q)^{1/2}} =(\sqrt{q!})^{d-1} \prod_{j=1}^d \frac{(H_j(2H_j-1))^{1/2}}{(q!(H'_j(2H'_j -1))^q)^{1/2}} 
\end{equation}
is the unique positive constant ensuring that $E[Z^{q, \HH}(\textbf{1})^2] = 1$, where
\begin{equation}\label{eq:HH}
\HH' := 1+ \frac{\HH-1}{q}\quad\big(\Longleftrightarrow (2\HH' -2)q = 2\HH-2\big).
\end{equation}
For a proof of (\ref{Def2}) when $d=2$, we refer to Tudor \cite[Chapter 4]{Tudor}. Extension to any value of $d$ as presented here is straightforward.

\section{Proof of Theorem \ref{main}}

We are now in a position to give the proof of our Theorem \ref{main}.
It is divided into three steps.

\subsection{Expanding into Wiener chaos}

In preparation of analysing the quadratic variation (\ref{eq:1}), let us find an explicit expression for the chaos decomposition of $V_\NN$. Using (\ref{Def2}) and proceeding by induction on the dimension $d$, we can write $\Delta Z^{q, \HH}_{[\frac{\ii}{\NN}, \frac{\ii +1}{\NN}]}$ as a $q$-th Wiener It\^{o} integral with respect to the standard two-sided Brownian sheet $(W_\t)_{\t \in \mathbb{R}^d}$ as follows: for every $0 \leq \ii \leq \NN-1$, one has
\begin{equation}\label{eq:Wchaos}
\Delta Z^{q, \HH}_{[\frac{\ii}{\NN}, \frac{\ii +1}{\NN}]} = I_q(f_{\ii,\NN}),
\end{equation} 
where
\begin{equation}\label{eq:Wchaos1}
f_{\ii,\NN}(\xx_1, \ldots, \xx_q) = b_{q, \HH}\prod_{j=1}^d f_{i_j,N_j}(x_{1,j}, \ldots, x_{q,j}),
\end{equation}
with $f_{i,N}(x_{1}, \ldots, x_{q})$ denoting the expression
\begin{align}\label{expression}
&\mathbbm{1}_{[0, \frac{i+1}{N}] }(x_{1} \vee \ldots \vee x_{q}) \int_{x_{1} \vee \ldots \vee x_{q}}^{\frac{i+1}{N}}du \partial_1K^{H'}(u, x_{1}) \ldots \partial_1K^{H'}(u, x_{q}) \nonumber\\
& - \mathbbm{1}_{[0, \frac{i}{N}] }(x_{1} \vee \ldots \vee x_{q}) \int_{x_{1} \vee \ldots \vee x_{q}}^{\frac{i}{N}}du \partial_1K^{H'}(u, x_{1}) \ldots \partial_1K^{H'}(u, x_{q}),
\end{align}
and with $b_{q, \HH}$ and $\HH'$ given by (\ref{eq:bqH}) and (\ref{eq:HH}) respectively. Indeed, for $d=1$, see \cite[Section 3, p.8]{TudorViens1}, it reduces to
$$\Delta Z^{q, H}_{[\frac{i}{N}, \frac{i +1}{N}]} = Z^{q, H}_{\frac{i+1}{N}} - Z^{q, H}_{\frac{i}{N}} =b_{q, H} I_q(f_{i,N}),$$
while for $d=2$, it is easy to verify that
$$\Delta Z^{q, \HH}_{[\frac{\ii}{\NN}, \frac{\ii +1}{\NN}]} = Z^{q, H_1, H_2}_{\frac{i+1}{N}, \frac{j+1}{M}} - Z^{q, H_1, H_2}_{\frac{i}{N}, \frac{j+1}{M}} - Z^{q, H_1, H_2}_{\frac{i+1}{N}, \frac{j}{M}} + Z^{q, H_1, H_2}_{\frac{i}{N}, \frac{j}{M}} = I_q(f_{i,j,N,M})$$
where
\begin{align*}
&f_{i,j,N,M} (x_1, y_1, \ldots, x_q, y_q)\\ 
& = b_{q, H_1, H_2} \mathbbm{1}_{[0, \frac{i+1}{N}] }(x_1 \vee \ldots \vee x_q) \int_{x_1 \vee \ldots \vee x_q}^{\frac{i+1}{N}}du \partial_1K^{H_1'}(u, x_1) \ldots \partial_1K^{H_1'}(u, x_q) \\
& \qquad\qquad\times \mathbbm{1}_{[0, \frac{j+1}{M}] }( y_1 \vee \ldots \vee y_q) \int_{y_1 \vee \ldots \vee y_q}^{\frac{j+1}{M}} dv \partial_1K^{H_2'}(v, y_1) \ldots \partial_1K^{H_2'}(v, y_q)\\
& - b_{q, H_1, H_2} \mathbbm{1}_{[0, \frac{i+1}{N}] }(x_1 \vee \ldots \vee x_q) \int_{x_1 \vee \ldots \vee x_q}^{\frac{i+1}{N}}du \partial_1K^{H_1'}(u, x_1) \ldots \partial_1K^{H_1'}(u, x_q) \\
& \qquad\qquad\times \mathbbm{1}_{[0, \frac{j}{M}] }( y_1 \vee \ldots \vee y_q) \int_{y_1 \vee \ldots \vee y_q}^{\frac{j}{M}} dv \partial_1K^{H_2'}(v, y_1) \ldots \partial_1K^{H_2'}(v, y_q)\\
& - b_{q, H_1, H_2} \mathbbm{1}_{[0, \frac{i}{N}] }(x_1 \vee \ldots \vee x_q) \int_{x_1 \vee \ldots \vee x_q}^{\frac{i}{N}}du \partial_1K^{H_1'}(u, x_1) \ldots \partial_1K^{H_1'}(u, x_q) \\
& \qquad\qquad\times \mathbbm{1}_{[0, \frac{j+1}{M}] }( y_1 \vee \ldots \vee y_q) \int_{y_1 \vee \ldots \vee y_q}^{\frac{j+1}{M}} dv \partial_1K^{H_2'}(v, y_1) \ldots \partial_1K^{H_2'}(v, y_q)\\
& + b_{q, H_1, H_2} \mathbbm{1}_{[0, \frac{i}{N}] }(x_1 \vee \ldots \vee x_q) \int_{x_1 \vee \ldots \vee x_q}^{\frac{i}{N}}du \partial_1K^{H_1'}(u, x_1) \ldots \partial_1K^{H_1'}(u, x_q) \\
& \qquad\qquad\times \mathbbm{1}_{[0, \frac{j}{M}] }( y_1 \vee \ldots \vee y_q) \int_{y_1 \vee \ldots \vee y_q}^{\frac{j}{M}} dv \partial_1K^{H_2'}(v, y_1) \ldots \partial_1K^{H_2'}(v, y_q)\\
& =  b_{q, H_1, H_2}f_{i,N}(x_1, \ldots, x_q)f_{j, M}(y_1, \ldots, y_q).
\end{align*}
The last equality above is obtained by grouping each term of $f_{i, j, N, M}$ together. Suppose that the expression (\ref{eq:Wchaos}), (\ref{eq:Wchaos1}) is true for $d$, that is, the kernel of $\Delta Z^{q, \HH}_{[\frac{\ii}{\NN}, \frac{\ii +1}{\NN}]}$ is equal to 
\eject
\begin{align*}
&b_{q, \HH} \sum_{(r_1, \ldots, r_d) \in \{0, 1\}^d} (-1)^{d - \sum_{i=1}^{d}r_i} \prod_{j=1}^d \mathbbm{1}_{[0, \frac{i_j + r_j}{N_j}] }(x_{1, j} \vee \ldots \vee x_{q, j}) \\
&\hspace{2cm} \times \int_{x_{1, j} \vee \ldots \vee x_{q, j}}^{\frac{i_j+r_j}{N_j}}du \partial_1K^{H_j'}(u, x_{1, j}) \ldots \partial_1K^{H_j'}(u, x_{q, j}) \\
&= b_{q, \HH} \prod_{j=1}^d f_{i_j, N_j}(x_{1, j}, \ldots, x_{q, j}).
\end{align*}
Then, for the case $d+1$ we have 
\begin{align*}
\Delta Z^{q, \HH}_{[\frac{\ii}{\NN}, \frac{\ii +1}{\NN}]}& = \sum_{\rr \in \{0, 1\}^{d+1}} (-1)^{d+1 - \sum_{i=1}^{d+1}r_i}Z^{q, \HH}_{\frac{\ii+\rr}{\NN}}\\ 
& = \sum_{(r_1, \ldots, r_d) \in \{0, 1\}^d} (-1)^{d - \sum_{i=1}^{d}r_i} Z^{q, \HH}_{\big(\frac{i_1+r_1}{N_1}, \ldots, \frac{i_d+r_d}{N_d}, \frac{i_{d+1}+1}{N_{d+1}} \big)} \\
& \qquad\qquad +  \sum_{(r_1, \ldots, r_d) \in \{0, 1\}^d} (-1)^{d +1 - \sum_{i=1}^{d}r_i} Z^{q, \HH}_{\big(\frac{i_1+r_1}{N_1}, \ldots, \frac{i_d+r_d}{N_d}, \frac{i_{d+1}}{N_{d+1}} \big)}\\
& = \sum_{(r_1, \ldots, r_d) \in \{0, 1\}^d} (-1)^{d - \sum_{i=1}^{d}r_i} \Big( Z^{q, \HH}_{\big(\frac{i_1+r_1}{N_1}, \ldots, \frac{i_d+r_d}{N_d}, \frac{i_{d+1}+1}{N_{d+1}} \big)} - Z^{q, \HH}_{\big(\frac{i_1+r_1}{N_1}, \ldots, \frac{i_d+r_d}{N_d}, \frac{i_{d+1}}{N_{d+1}} \big)}  \Big).
\end{align*}
It belongs to the $q$-Wiener chaos with the kernel  $f_{\ii, \NN}$ given by 
\begin{align*}
f_{\ii, \NN}& = b_{q, \HH} \sum_{(r_1, \ldots, r_d) \in \{0, 1\}^d} (-1)^{d - \sum_{i=1}^{d}r_i} \prod_{j=1}^d \mathbbm{1}_{[0, \frac{i_j + r_j}{N_j}] }(x_{1, j} \vee \ldots \vee x_{q, j})\\
& \hspace{5cm} \times \int_{x_{1, j} \vee \ldots \vee x_{q, j}}^{\frac{i_j+r_j}{N_j}}du \partial_1K^{H_j'}(u, x_{1, j}) \ldots \partial_1K^{H_j'}(u, x_{q, j}) \\ 
& \times \bigg(\int_{x_{1, d+1} \vee \ldots \vee x_{q, d+1}}^{\frac{i_{d+1}+ 1}{N_{d+1}}}du' \partial_1K^{H'_{d+1}}(u', x_{1, d+1}) \ldots \partial_1K^{H'_{d+1}}(u', x_{q, d+1}) \\
& \hspace{2.5cm} - \int_{x_{1, d+1} \vee \ldots \vee x_{q, d+1}}^{\frac{i_{d+1}}{N_{d+1}}}du' \partial_1K^{H'_{d+1}}(u', x_{1, d+1}) \ldots \partial_1K^{H'_{d+1}}(u', x_{q, d+1})\bigg).
\end{align*}
By the induction hypothesis, one gets
$f_{\ii, \NN} =  b_{q, \HH}  \prod_{j=1}^{d+1} f_{i_j, N_j}(x_{1, j}, \ldots, x_{q, j}),$
which is our desired expression.

Next, by applying the product formula (\ref{eq:P1}), we  can write
\begin{equation}\label{eq:dola1}
\Big(\Delta Z^{q, \HH}_{[\frac{\ii}{\NN}, \frac{\ii +1}{\NN}]}\Big)^2 - E\Big[\Big(\Delta Z^{q, \HH}_{[\frac{\ii}{\NN}, \frac{\ii +1}{\NN}]}\Big)^2 \Big]= \sum_{r=0}^{q-1} r!\binom{q}{r}^2 I_{2q-2r}(f_{\ii,\NN} \widetilde{\otimes}_r f_{\ii,\NN}).
\end{equation}
Let us compute the contractions appearing in the right-hand side of (\ref{eq:dola1}). For every $0 \leq r \leq q-1$, we have
{\allowdisplaybreaks
\begin{align}\label{eq:3.1.1}
&(f_{\ii,\NN} \otimes_r f_{\ii,\NN})(\xx_1, \ldots, \xx_{2q-2r}) \nonumber \\ 
& = \int_{([0, 1]^{d})^{  r}}d\aa_1\ldots d\aa_r f_{\ii,\NN}(\xx_1, \ldots, \xx_{q-r}, \aa_1, \ldots, \aa_r) \nonumber\\
&\qquad\qquad\qquad\qquad \times f_{\ii,\NN}(\xx_{q-r+1}, \ldots, \xx_{2q-2r}, \aa_1, \ldots, \aa_r) \nonumber\\
&= b_{q, \HH}^2 \int_{([0, 1]^{d })^{ r}}d\aa_1\ldots d\aa_r \prod_{j=1}^d f_{i_j, N_j}(x_{1,j}, \ldots, x_{q-r,j}, a_{1,j}, \ldots, a_{r,j}) \nonumber\\
&\qquad\qquad\qquad\qquad \times  \prod_{j=1}^d f_{i_j, N_j}(x_{q-r+1,j}, \ldots, x_{2q-2r,j}, a_{1,j}, \ldots, a_{r,j}) \nonumber\\
& = b_{q, \HH}^2 \prod_{j=1}^d (f_{i_j ,N_j} \otimes_r f_{i_j ,N_j})(x_{1,j}, \ldots, x_{2q-2r,j}),
\end{align}
}where 
{\allowdisplaybreaks
\begin{align}\label{eq:3.1}
&(f_{i, N} \otimes_r f_{i, N})(x_{1}, \ldots, x_{2q-2r})= (H'(2H' -1))^r \nonumber\\
& \times \bigg\{\mathbbm{1}_{[0, \frac{i_+1}{N}] }(x_{1}\vee \ldots x_{q-r}) \int_{x_{1}\vee \ldots x_{q-r}}^{\frac{i+1}{N}}du \partial_1K^{H'}(u, x_{1})\ldots \partial_1K^{H'}(u, x_{q-r}) \nonumber\\
&\qquad\times\mathbbm{1}_{[0, \frac{i+1}{N}]}(x_{q-r+1}\vee \ldots x_{2q-2r}) \int_{x_{q-r+1}\vee \ldots x_{2q-2r}}^{\frac{i+1}{N}}du' \partial_1K^{H'}(u', x_{q-r+1})\ldots \nonumber \\
&\hspace{8cm}\ldots \partial_1K^{H'}(u', x_{2q-2r}) |u-u'|^{(2H'-2)r} \nonumber\\
&\quad -\mathbbm{1}_{[0, \frac{i+1}{N}] }(x_{1}\vee \ldots x_{q-r}) \int_{x_{1}\vee \ldots x_{q-r}}^{\frac{i+1}{N}}du \partial_1K^{H'}(u, x_{1})\ldots \partial_1K^{H'}(u, x_{q-r})\nonumber\\
&\qquad\times\mathbbm{1}_{[0, \frac{i}{N}]}(x_{q-r+1}\vee \ldots x_{2q-2r}) \int_{x_{q-r+1}\vee \ldots x_{2q-2r}}^{\frac{i}{N}}du' \partial_1K^{H'}(u', x_{q-r+1})\ldots  \nonumber\\
&\hspace{8cm} \ldots \partial_1K^{H'}(u', x_{2q-2r}) |u-u'|^{(2H'-2)r} \nonumber\\
&\quad -  \mathbbm{1}_{[0, \frac{i}{N}] }(x_{1}\vee \ldots x_{q-r}) \int_{x_{1}\vee \ldots x_{q-r}}^{\frac{i}{N}}du \partial_1K^{H'}(u, x_{1})\ldots \partial_1K^{H'}(u, x_{q-r}) \nonumber\\
&\qquad\times\mathbbm{1}_{[0, \frac{i+1}{N}]}(x_{q-r+1}\vee \ldots x_{2q-2r}) \int_{x_{q-r+1}\vee \ldots x_{2q-2r}}^{\frac{i+1}{N}}du' \partial_1K^{H'}(u', x_{q-r+1})\ldots \nonumber \\
&\hspace{8cm}\ldots \partial_1K^{H'}(u', x_{2q-2r}) |u-u'|^{(2H'-2)r} \nonumber\\
&\quad+\mathbbm{1}_{[0, \frac{i}{N}] }(x_{1}\vee \ldots x_{q-r}) \int_{x_{1}\vee \ldots x_{q-r}}^{\frac{i}{N}}du \partial_1K^{H'}(u, x_{1})\ldots \partial_1K^{H'}(u, x_{q-r}) \nonumber\\
&\qquad\times\mathbbm{1}_{[0, \frac{i}{N}]}(x_{q-r+1}\vee \ldots x_{2q-2r}) \int_{x_{q-r+1}\vee \ldots x_{2q-2r}}^{\frac{i}{N}}du' \partial_1K^{H'}(u', x_{q-r+1})\ldots \nonumber\\
&\hspace{8cm}\ldots \partial_1K^{H'}(u', x_{2q-2r}) |u-u'|^{(2H'-2)r}\bigg\}.
\end{align}
}(See \cite[page 10]{TudorViens1} for a detailed computation of the expression (\ref{eq:3.1}).) Moreover, since $Z^{q,\HH}$ is $\HH$-self-similar and has stationary increments, one has
$$\Delta Z^{q, \HH}_{[\frac{\ii}{\NN}, \frac{\ii +1}{\NN}]} \overset{(d)}{=} \NN^{-\HH} \Delta Z^{q, \HH}_{[\ii, \ii+1]}\overset{(d)}{=} \NN^{-\HH}Z^{q, \HH}_{[0,\bf{1}]}.$$
It follows that 
$$E\bigg[\NN^{2\HH}\Big(\Delta Z^{q, \HH}_{[\frac{\ii}{\NN}, \frac{\ii +1}{\NN}]}\Big)^2\bigg] = E[Z^{q, \HH}(\textbf{1})^2] = 1.$$
As a consequence, we have 
\begin{equation}\label{eq:VN}
V_{\NN} = F_{2q, \NN} + c_{2q-2}F_{2q-2, \NN} + \ldots + c_4F_{4, \NN} + c_2F_{2, \NN}.
\end{equation}
where
$c_{2q-2r} = r! \binom{q}{r}^2$, $r=0,\ldots, q-1$,
are the combinator constants coming from the product formula, and
\begin{equation}\label{eq:3.4}
F_{2q-2r, \NN}: = \NN^{2\HH -1}I_{2q-2r}\bigg(\sum_{\ii=0}^{\NN-1} f_{\ii,\NN} \widetilde{\otimes}_r f_{\ii,\NN}\bigg),
\end{equation}
for the kernels $f_{\ii,\NN} \otimes_r f_{\ii,\NN}$ computed in (\ref{eq:3.1.1})-(\ref{eq:3.1}).

\subsection{Evaluating the $L^2(\Omega)$-norm}

Set
\begin{equation}\label{eq:3.5}
c_{1, \HH}= \frac{2! 2^d b_{q, \HH}^4 (\HH'(2\HH'-1))^{2q}}{(4\HH' - 3)(4\HH'-2)[(2\HH'-2)(q-1)+1]^2[(\HH'-1)(q-1) + 1]^2}.
\end{equation}
We claim that 
\begin{equation}\label{eq:30}
\lim_{\NN \to \infty} E[c_{1, \HH}^{-1}\NN^{2(2-2\HH')}c_2^{-2}V_{\NN}^2] = 1.
\end{equation}

Let us prove (\ref{eq:30}). Due to the orthogonality property for Wiener chaoses of different orders, it is sufficient to evaluate the $L^2(\Omega)$-norm of each multiple Wiener-It\^{o} integrals appearing in the chaotic decomposition (\ref{eq:VN}) of $V_\NN$. Let us start with the double integral:
$$F_{2, \NN}= \NN^{2\HH-1}I_2\bigg(\sum_{\ii=0}^{\NN-1} f_{\ii,\NN} \otimes_{q-1} f_{\ii,\NN}\bigg).$$
Since the kernel $\sum_{\ii=0}^{\NN-1} f_{\ii,\NN} \otimes_{q-1} f_{\ii,\NN}$ is symmetric, one has
\begin{align*}
E[F_{2, \NN}^2]& = 2!N^{4\HH-2} \bigg\| \sum_{\ii=0}^{\NN-1} f_{\ii,\NN} \otimes_{q-1} f_{\ii,\NN}\bigg\|_{L^2(([0, 1]^{d})^{ 2})}^2\\ 
&=2!\NN^{4\HH-2} \sum_{\ii, \kk=0}^{\NN-1}\left\langle  f_{\ii,\NN} \otimes_{q-1} f_{\ii,\NN},  f_{\kk,\NN} \otimes_{q-1} f_{\kk,\NN} \right\rangle_{L^2(([0, 1]^{d})^{ 2})}.
\end{align*}
Let us now compute the scalar products in the above expression. By using (\ref{eq:3.1.1}), (\ref{eq:3.1}), by applying Fubini's theorem and by noting that $\int_0^{u \wedge v} \partial_1 K^{H'}(u, a) \partial_1 K^{H'}(v, a)da = H'(2H' -1)|u-v|^{2H'-2}$, it is easy to verify that
\begin{align*}
&\left\langle  f_{\ii,\NN} \otimes_{q-1} f_{\ii,\NN},  f_{\kk,\NN} \otimes_{q-1} f_{\kk,\NN} \right\rangle_{L^2(([0, 1]^{d})^{ 2})}\\ 
&=b_{q, \HH}^4 \prod_{j=1}^d \left\langle f_{i_j,N_j} \otimes_{q-1} f_{i_j,N_j}, f_{k_j,N_j} \otimes_{q-1} f_{k_j,N_j} \right\rangle_{L^2([0, 1]^2)}\\
& = b^4_{q, \HH} (\HH'(2\HH'-1))^{2q} \prod_{j=1}^d \int_{\frac{i_j}{N_j}}^{\frac{i_j+1}{N_j}}du_j\int_{\frac{i_j}{N_j}}^{\frac{i_j+1}{N_j}}dv_j \int_{\frac{k_j}{N_j}}^{\frac{k_j+1}{N_j}}du'_j \int_{\frac{k_j}{N_j}}^{\frac{k_j+1}{N_j}}dv'_j \\
& \hspace{4.5cm}\times |u_j-v_j|^{(2H'_j-2)(q-1)}|u'_j-v'_j|^{(2H'_j-2)(q-1)}\\
&\hspace{4.5cm}\times|u_j-u'_j|^{2H'_j-2}|v_j-v'_j|^{2H'_j-2},
\end{align*}
(see, e.g., \cite[page 11]{TudorViens1}). The change of variables $u' = (u-\frac{i}{N})N$ for each $u_j, u'_j, v_j, v'_j$ with $j$ from $1$ to $d$ yields
\begin{align}\label{f1}
E[F_{2, \NN}^2] &= 2 b^4_{q, \HH} (\HH'(2\HH'-1))^{2q}\NN^{4\HH-2}\NN^{-4}\NN^{-(2\HH'-2)2q} \nonumber\\ 
&\times  \sum_{\ii, \kk=0}^{\NN-1}\prod_{j=1}^d \int_0^1du_j\int_0^1dv_j \int_0^1du'_j \int_0^1dv'_j  
|u_j-v_j|^{(2H'_j-2)(q-1)}|u'_j-v'_j|^{(2H'_j-2)(q-1)} \nonumber\\
&\hspace{5cm}\times|u_j-u'_j + i_j -k_j|^{2H'_j-2}|v_j-v'_j + i_j -k_j|^{2H'_j-2}.
\end{align}
Now, we split the sum $\sum_{\ii, \kk=0}^{\NN-1}$ appearing in $E[F_{2, \NN}^2]$ just above into
\begin{equation}\label{sum}
\sum_{\ii, \kk=0}^{\NN-1} = \sum_{\substack{\ii, \kk= 0 \\ \exists  1\leq j \leq d: i_j = k_j}}^{\NN-1} + \sum_{\substack{\ii, \kk= 0 \\ \forall  j: i_j \ne k_j}}^{\NN-1} .
\end{equation}
For the first term in the right-hand side of (\ref{sum}), without loss of generality, let us assume that $ i_1 = k_1, \ldots, i_m = k_m$ for some $1 \leq m < d$ and $i_j \ne k_j$ for all $ m+1 \leq j \leq d$. Then,
\begin{align*}
&\NN^{-2}\sum_{\substack{\ii, \kk= 0 \\ i_1 = k_1, \ldots, i_m = k_m}}^{\NN-1}\prod_{j=1}^d \int_{[0, 1]^4}du_jdv_jdu'_jdv'_j |u_j-v_j|^{(2H'_j-2)(q-1)}|u'_j-v'_j|^{(2H'_j-2)(q-1)}\\
&\hspace{5.5cm}\times|u_j-u'_j + i_j -k_j|^{2H'_j-2}|v_j-v'_j + i_j -k_j|^{2H'_j-2}\\
&= \prod_{j=1}^m N_j^{-1}  \int_{[0,1]^4}du_jdv_jdu'_jdv'_j (|u_j-v_j||u'_j-v'_j|)^{(2H'_j-2)(q-1)}(|u_j-u'_j ||v_j-v'_j |)^{2H'_j-2}\\
& \qquad\times \sum_{\substack{i_{m+1}, k_{m+1} =0\\ i_{m+1} \ne k_{m+1}}}^{N_j-1} \ldots \sum_{\substack{i_d, k_d =0\\ i_d \ne k_d}}^{N_j-1} \prod_{j=m+1}^d \int_{[0,1]^4}du_jdv_jdu'_jdv'_j (|u_j-v_j||u'_j-v'_j|)^{(2H'_j-2)(q-1)}\\
& \hspace{5cm}\times N_j^{-2}|u_j-u'_j + i_j -k_j|^{2H'_j-2}|v_j-v'_j + i_j -k_j|^{2H'_j-2}.
\end{align*}
By switching sum and product in the above expression, we arrive 
{\allowdisplaybreaks
\begin{align*}
& \prod_{j=1}^m N_j^{-1}  \int_{[0,1]^4}du_jdv_jdu'_jdv'_j (|u_j-v_j||u'_j-v'_j|)^{(2H'_j-2)(q-1)}(|u_j-u'_j ||v_j-v'_j |)^{2H'_j-2}\\
& \qquad\times \prod_{j= m+1}^d \bigg(\sum_{\substack{i_j, k_j = 0 \\i_j \ne k_j}}^{N_j-1} \int_{[0,1]^4}du_jdv_jdu'_jdv'_j |u_j-v_j|^{(2H'_j-2)(q-1)}|u'_j-v'_j|^{(2H'_j-2)(q-1)}\\
& \qquad\qquad\qquad\times N_j^{-2}|u_j-u'_j + i_j -k_j|^{2H'_j-2}|v_j-v'_j + i_j -k_j|^{2H'_j-2}\bigg)\\
&=  \prod_{j=1}^m N_j^{-1}  \int_{[0,1]^4}du_jdv_jdu'_jdv'_j (|u_j-v_j||u'_j-v'_j|)^{(2H'_j-2)(q-1)}(|u_j-u'_j ||v_j-v'_j |)^{2H'_j-2}\\
& \qquad\times \prod_{j=m+1}^d \bigg( \int_{[0,1]^4}du_jdv_jdu'_jdv'_j |u_j-v_j|^{(2H'_j-2)(q-1)}|u'_j-v'_j|^{(2H'_j-2)(q-1)}\\
& \qquad\qquad\qquad\times 2N_j^{-2}\sum_{\substack{i_j, k_j= 0 \\ i_j > k_j}}^{N_j-1}|u_j-u'_j + i_j -k_j|^{2H'_j-2}|v_j-v'_j + i_j -k_j|^{2H'_j-2}\bigg).
\end{align*}
}One has that
\begin{align*}
&N^{-2}\sum_{\substack{i, k= 0 \\ i > k}}^{N-1}|u-u' + i -k|^{2H'-2}|v-v' + i -k|^{2H'-2}\\ 
&= N^{2(2H'-2)}\frac{1}{N}\sum_{n=1}^{N}\Big(1-\frac{n}{N}\Big)\Big|\frac{u-u'}{N}+ \frac{n}{N}\Big|^{2H'-2}\Big|\frac{v-v'}{N}+\frac{n}{N}\Big|^{2H'-2}
\end{align*}
 is asymptotically equivalent to 
$N^{2(2H'-2)}\int_0^1(1-x)x^{4H'-4}dx = N^{2(2H'-2)}\frac{1}{(4H' - 3)(4H'-2)}.$
It follows that
\begin{align*}
&\NN^{-2}\sum_{\substack{\ii, \kk= 0 \\ i_1 = k_1, \ldots, i_m = k_m}}^{\NN-1}\prod_{j=1}^d \int_{[0, 1]^4}du_jdv_jdu'_jdv'_j |u_j-v_j|^{(2H'_j-2)(q-1)}|u'_j-v'_j|^{(2H'_j-2)(q-1)}\\
&\hspace{4.5cm}\times|u_j-u'_j + i_j -k_j|^{2H'_j-2}|v_j-v'_j + i_j -k_j|^{2H'_j-2}\\
 &\approx \prod_{j=1}^m N_j^{-1}  \int_{[0,1]^4}du_jdv_jdu'_jdv'_j |u_j-v_j|^{(2H'_j-2)(q-1)}|u'_j-v'_j|^{(2H'_j-2)(q-1)}\\
&\hspace{5.5cm}\times|u_j-u'_j |^{2H'_j-2}|v_j-v'_j |^{2H'_j-2}\\
& \times \prod_{j=m+1}^d 2N_j^{2(2H'_j-2)}\frac{1}{(4H'_j - 3)(4H'_j-2)} \bigg( \int_{[0,1]^2}du_jdv_j |u_j-v_j|^{(2H'_j-2)(q-1)}\bigg)^2.
\end{align*}
Since $2(2-2H'_j) - 1 < 0$ for all $j$, one gets, as $\NN \to \infty$,
\begin{align}\label{f2}
&\NN^{2(2-2\HH'_j)} \times \NN^{-2}\sum_{\substack{\ii, \kk= 0 \\ \exists  1\leq j \leq d: i_j = k_j}}^{\NN-1}\prod_{j=1}^d \int_{[0, 1]^4}du_jdv_jdu'_jdv'_j |u_j-v_j|^{(2H'_j-2)(q-1)}|u'_j-v'_j|^{(2H'_j-2)(q-1)}\nonumber\\
&\hspace{4.5cm}\times |u_j-u'_j + i_j -k_j|^{2H'_j-2}|v_j-v'_j + i_j -k_j|^{2H'_j-2} \longrightarrow  0.
 \end{align}
Similarly for the second term in (\ref{sum}), that is, when $i_j \ne k_j$ for all $1 \leq j \leq d$, we have
\begin{align*}
&\NN^{-2}\sum_{\substack{\ii, \kk= 0 \\ i_j \ne k_j, \text{ } \forall j}}^{\NN-1}\prod_{j=1}^d \int_{[0, 1]^4}du_jdv_jdu'_jdv'_j |u_j-v_j|^{(2H'_j-2)(q-1)}|u'_j-v'_j|^{(2H'_j-2)(q-1)}\\
&\hspace{4.5cm}\times|u_j-u'_j + i_j -k_j|^{2H'_j-2}|v_j-v'_j + i_j -k_j|^{2H'_j-2}\\
 &\approx \prod_{j=1}^d  N_j^{2(2H'_j-2)}\frac{2}{(4H'_j - 3)(4H'_j-2)}\bigg( \int_{[0,1]^2}du_jdv_j |u_j-v_j|^{(2H'_j-2)(q-1)}\bigg)^2\\
 &= \prod_{j=1}^d N_j^{2(2H'_j-2)}\frac{2}{(4H'_j - 3)(4H'_j-2)[(2H'_j-2)(q-1)+1]^2[(H'_j-1)(q-1) + 1]^2}.
\end{align*}
It follows that
\begin{align}\label{f3}
&\NN^{2(2-2\HH')} \times \NN^{-2}\sum_{\substack{\ii, \kk= 0 \\ i_j \ne k_j, \text{ } \forall j}}^{\NN-1}\prod_{j=1}^d \int_{[0, 1]^4}du_jdv_jdu'_jdv'_j |u_j-v_j|^{(2H'_j-2)(q-1)}|u'_j-v'_j|^{(2H'_j-2)(q-1)}\nonumber\\
&\hspace{5.5cm}\times|u_j-u'_j + i_j -k_j|^{2H'_j-2}|v_j-v'_j + i_j -k_j|^{2H'_j-2}\nonumber\\
 &\longrightarrow \prod_{j=1}^d  \frac{2}{(4H'_j - 3)(4H'_j-2)[(2H'_j-2)(q-1)+1]^2[(H'_j-1)(q-1) + 1]^2}.
\end{align}

To conclude that 
\begin{equation}\label{eq:3.6}
\lim_{\NN \to \infty} E[c_{1, \HH}^{-1}\NN^{2(2-2\HH')}F_{2, \NN}^2] = 1,
\end{equation}
we use the expression  (\ref{f1}) for $E[F_{2,\NN}^2]$.
The first sum in (\ref{sum}) goes to zero according to (\ref{f2}), whereas the second sum goes to
the quantity in (\ref{f3}).
Going back to the definition (\ref{eq:3.5}) of $c_{1,{\bf H}}$, we arrive to the desired conclusion (\ref{eq:3.6}).

Let us now consider the remaining terms $F_{4, \NN}, \ldots, F_{2q, \NN}$ in the chaos decomposition (\ref{eq:VN}). 
Using that $\|\widetilde{g}\|_{L^2} \leq \|g\|_{L^2}$ for any square integrable function $g$, one can write, for every $ 0 \leq r \leq q-2$,
\begin{align*}
E[F_{2q-2r, \NN}^2]& = \NN^{4\HH-2}(2q-2r)! \,\bigg\| \sum_{\ii=0}^{\NN-1} f_{\ii,\NN} \widetilde{\otimes}_r f_{\ii,\NN}\bigg\|^2_{L^2(([0, 1]^{d})^{2q-2r})}\\ 
 & \leq \NN^{4\HH-2}(2q-2r)!\,\bigg\| \sum_{\ii=0}^{\NN-1}f_{\ii,\NN} \otimes_r f_{\ii,\NN}\bigg\|^2_{L^2(([0, 1]^{d})^{2q-2r})}\\ 
&= (2q-2r)!\NN^{4\HH-2} \sum_{\ii, \kk=0}^{\NN-1}\left\langle  f_{\ii,\NN} \otimes_r f_{\ii,\NN},  f_{\kk,\NN} \otimes_r f_{\kk,\NN} \right\rangle_{L^2(([0, 1]^{d})^{2q-2r})}.
\end{align*}
Proceeding as above, we obtain 
\begin{align*}
&\left\langle  f_{\ii,\NN} \otimes_r f_{\ii,\NN},  f_{\kk,\NN} \otimes_r f_{\kk,\NN} \right\rangle_{L^2([0, 1]^{d \cdot (2q-2r)})}\\ 
& = b^4_{q, \HH} (\HH'(2\HH'-1))^{2q} \prod_{j=1}^d \int_{\frac{i_j}{N_j}}^{\frac{i_j+1}{N_j}}du_j\int_{\frac{i_j}{N_j}}^{\frac{i_j+1}{N_j}}dv_j \int_{\frac{k_j}{N_j}}^{\frac{k_j+1}{N_j}}du'_j \int_{\frac{k_j}{N_j}}^{\frac{k_j+1}{N_j}}dv'_j \\
& \hspace{4.5cm}\times |u_j-v_j|^{(2H'_j-2)r}|u'_j-v'_j|^{(2H'_j-2)r}\\
&\hspace{4.5cm}\times|u_j-u'_j|^{(2H'_j-2)(q-r)}|v_j-v'_j|^{(2H'_j-2)(q-r)}.
\end{align*}
Using the change of variables $u' = (u-\frac{i}{N})N$ for each $u_j, u_j, v_j, v'_j$ with $j=1, \ldots, d$, one obtains
\begin{align*}
E[F_{2q-2r, \NN}^2] &\leq (2q-2r)! b^4_{q, \HH} (\HH'(2\HH'-1))^{2q}\NN^{4\HH-2}\NN^{-4}\NN^{-(2\HH'-2)2q}\\ 
&\times  \sum_{\ii, \kk=0}^{\NN-1}\bigg(\prod_{j=1}^d \int_0^1du_j\int_0^1dv_j \int_0^1du'_j \int_0^1dv'_j |u_j-v_j|^{(2H'_j-2)r}|u'_j-v'_j|^{(2H'_j-2)r}\\
&\hspace{2cm}\times|u_j-u'_j + i_j -k_j|^{(2H'_j-2)(q-r)}|v_j-v'_j + i_j -k_j|^{(2H'_j-2)(q-r)}\bigg).
\end{align*}
Switching sum and product in the above expression, one obtains
\begin{align}\label{Khochiu}
E[F_{2q-2r, \NN}^2] &\leq (2q-2r)! b^4_{q, \HH} (\HH'(2\HH' -1))^{2q}\NN^{-2} \nonumber\\
&\qquad\times\prod_{j=1}^d\int_{[0, 1]^4}du_jdv_jdu'_jdv'_j |u_j-v_j|^{(2H'_j-2)r}|u'_j-v'_j|^{(2H'_j-2)r}\nonumber\\ 
&\qquad\times \bigg( \sum_{i_j, k_j=0}^{N_j-1}|u_j-u'_j + i_j -k_j|^{(2H'_j-2)(q-r)}|v_j-v'_j + i_j -k_j|^{(2H'_j-2)(q-r)}\bigg).
\end{align}
Note that the above sum $\sum_{i_j, k_j=0}^{N_j-1}$ can be divided into two parts: the diagonal part with $i_j=k_j$ and the non-diagonal part with $i_j\ne k_j$. It is easily seen that the non-diagonal part is dominant. Indeed, the diagonal part in the right-hand side of (\ref{Khochiu}) is equal to
\begin{align*}
&(2q-2r)! b^4_{q, \HH} (\HH'(2\HH' -1))^{2q}\NN^{-1}\prod_{j=1}^d\int_{[0, 1]^4}du_jdv_jdu'_jdv'_j\\
&\qquad\times |u_j-v_j|^{(2H'_j-2)r}|u'_j-v'_j|^{(2H'_j-2)r}|u_j-u'_j |^{(2H'_j-2)(q-r)}|v_j-v'_j|^{(2H'_j-2)(q-r)}.
\end{align*}
and it tends to zero since $(2H'_j-2)r >- 1$ and $(2H'_j-2)(q-r) > -1$. Thus, in order to find a bound of $E[F_{2q-2r, \NN}^2]$ in (\ref{Khochiu}), we have to study the following sum
\begin{equation}\label{Khochiu1}
\frac{1}{N^2}\sum_{\substack{i, k=0 \\ i \ne k}}^{N-1}|u-u' + i -k|^{(2H'-2)(q-r)}|v-v' + i -k|^{(2H'-2)(q-r)}
\end{equation}
for all $q \geq 2$ and $r =0, \ldots, q-2$, when $u, u', v, v' \in [0, 1]$. In (\ref{Khochiu1}), one has set $H' = 1 + \frac{H-1}{q}$ with $H > \frac12$. We now analyse the behavior of (\ref{Khochiu1}) according to the following three cases: $H > \frac34, H <\frac34$ and $H =\frac34$.

$\bullet$ If $H > \frac34$, then (\ref{Khochiu1}) is equal to
$$ N^{(2H'-2)(2q-2r)}\frac{2}{N}\sum_{n=1}^{N}\Big(1-\frac{n}{N}\Big)\Big|\frac{u-u'}{N}+ \frac{n}{N}\Big|^{(2H'-2)(q-r)}\Big|\frac{v-v'}{N}+\frac{n}{N}\Big|^{(2H'-2)(q-r)}.$$
By multiplying (\ref{Khochiu1}) by $N^{(2-2H')(2q-2r)}$ one has
\begin{align*}
&N^{(2-2H')(2q-2r)}\times \frac{1}{N^2}\sum_{\substack{i, k=0 \\ i \ne k}}^{N-1}|u-u' + i -k|^{(2H'-2)(q-r)}|v-v' + i -k|^{(2H'-2)(q-r)}\\ 
& = \frac{2}{N}\sum_{n=1}^{N}\Big(1-\frac{n}{N}\Big)\Big|\frac{u-u'}{N}+ \frac{n}{N}\Big|^{(2H'-2)(q-r)}\Big|\frac{v-v'}{N}+\frac{n}{N}\Big|^{(2H'-2)(q-r)}\\
&\approx 2 \int_0^1 (1-x)x^{2(2H' -2)(q-r)}dx < \infty \qquad \text{since } H  > \frac34.
\end{align*}

$\bullet$ If $H <\frac34$, (\ref{Khochiu1}) is bounded by
$$\frac{1}{N}\sum_{r\in\mathbb{Z}\setminus\{0\}} |u-u' + r|^{(2H'-2)(2q-2r)}|v-v' + r|^{(2H'-2)(2q-2r)} = O(\frac{1}{N}).$$

$\bullet$ If $H =\frac34$, following the same route as in the case $H < \frac34$, we arrive to $(\ref{Khochiu1}) = O(\frac{\log{\NN}}{\NN})$. 

Now, we go back to (\ref{Khochiu}). From the analysis of (\ref{Khochiu1}), we conclude that
$$E[F_{2q-2r, \NN}^2] = \begin{cases}
O(\NN^{-(2\HH' -2)(2q-2r)})& \text{if } H \in (\frac34, 1)\\
O(\NN^{-1})& \text{if } H \in (\frac12, \frac34)\\
O(\frac{\log{\NN}}{\NN})& \text{if } H =\frac34\\
\end{cases}
$$
Therefore, for all $0 \leq r \leq q-2$ and as $\NN \to \infty$, one has
\begin{equation}\label{eq:3.7}
\lim_{\NN \to \infty} E[\NN^{2(2-2\HH')}F_{2q-2r, \NN}^2] = 0.
\end{equation}
Thus, from (\ref{eq:3.6}), (\ref{eq:3.7}) and the orthogonality of Wiener chaos, we obtain (\ref{eq:30}).

\subsection{Concluding the proof of Theorem \ref{main}}

Thanks to (\ref{eq:3.7}), in order to understand the asymptotic behavior of  the normalized sequence of $V_{\NN}$, it is enough to analyse the convergence of the term

\begin{equation}\label{eq:3.12}
\NN^{2-2\HH'}F_{2, \NN} = I_2\bigg(\NN^{2\HH-1}\NN^{2-2\HH'}\sum_{\ii=0}^{\NN-1}f_{\ii,\NN} \otimes_{q-1} f_{\ii,\NN} \bigg),
\end{equation}
with
{\allowdisplaybreaks
\begin{align*}
&f_{\ii,\NN} \otimes_{q-1} f_{\ii,\NN} (\xx_1, \xx_2)= b_{q, \HH}^2 \prod_{j=1}^d(f_{i_j, N_j} \otimes_{q-1} f_{i_j, N_j})(x_{1,j}, x_{2,j})\\
& = b_{q, \HH}^2 (\HH'(2\HH'-1))^{q-1}\\
&\qquad \times \prod_{j=1}^d\bigg(\mathbbm{1}_{[0, \frac{i_j}{N_j}]}(x_{1,j}) \mathbbm{1}_{[0, \frac{i_j}{N_j}]}(x_{2,j}) \int_{\frac{i_j}{N_j}}^{\frac{i_j+1}{N_j}}du\int_{\frac{i_j}{N_j}}^{\frac{i_j+1}{N_j}}du'  \partial_1K^{H'_j}(u, x_{1,j}) \\
&\hspace{7cm} \times \partial_1K^{H'_j}(u', x_{2,j})|u-u'|^{(2H'_j-2)(q-1)}\\
&\hspace{1.2cm}+ \mathbbm{1}_{[0, \frac{i_j}{N_j}]}(x_{1,j}) \mathbbm{1}_{[\frac{i_j}{N_j}, \frac{i_j+1}{N_j}]}(x_{2,j}) \int_{\frac{i_j}{N_j}}^{\frac{i_j+1}{N_j}}du\int_{x_{2,j}}^{\frac{i_j+1}{N_j}}du'  \partial_1K^{H'_j}(u, x_{1,j})\\
&\hspace{7cm}\times \partial_1K^{H'_j}(u', x_{2,j})|u-u'|^{(2H'_j-2)(q-1)}\\
&\hspace{1.2cm}+\mathbbm{1}_{[\frac{i_j}{N_j}, \frac{i_j+1}{N_j}]}(x_{1,j}) \mathbbm{1}_{[0, \frac{i_j+1}{N_j}]}(x_{2,j}) \int_{x_{1,j}}^{\frac{i_j+1}{N_j}}du\int_{\frac{i_j}{N_j}}^{\frac{i_j+1}{N_j}}du'  \partial_1K^{H'_j}(u, x_{1,j})\\
&\hspace{7cm} \times\partial_1K^{H'_j}(u', x_{2,j})|u-u'|^{(2H'_j-2)(q-1)}\\
&\hspace{1.2cm}+\mathbbm{1}_{[\frac{i_j}{N_j}, \frac{i_j}{N_j}]}(x_{1,j}) \mathbbm{1}_{[\frac{i_j}{N_j}, \frac{i_j+1}{N_j}]}(x_{2,j}) \int_{x_{1,j}}^{\frac{i_j+1}{N_j}}du\int_{x_{2,j}}^{\frac{i_j+1}{N_j}}du'  \partial_1K^{H'_j}(u, x_{1,j}) \\
&\hspace{7cm} \times\partial_1K^{H'_j}(u', x_{2,j})|u-u'|^{(2H'_j-2)(q-1)}\bigg).
\end{align*}
}Among the four terms in the right-hand side of the above expression, only the first one is not asymptotically negligible in $L^2(\Omega)$ as $\NN \to \infty$, see \cite[page 14 and 15]{TudorViens1} or follow the lines of \cite{TudorViens2} for details. Furthermore, by the isometry property for multiple Wiener-It\^{o} integrals, evaluating the $L^2(\Omega)$-limit  of a sequence belonging to the second Wiener chaos is equivalent to evaluating the $L^2(([0, 1]^{d})^{2})$-limit of the sequence of their corresponding symmetric kernels. Therefore, we are left to find the limit of $f_2^{\NN}$ in $L^2(([0, 1]^{d})^2)$, where
\eject
\begin{align*}
f_2^{\NN}(\xx_1, \xx_2):&= \NN^{2\HH-1} \NN^{2-2\HH'}b_{q, \HH}^2 (\HH'(2\HH'-1))^{q-1}\\
& \times \sum_{\ii =0}^{\NN -1}\bigg(  \prod_{j=1}^d \mathbbm{1}_{[0, \frac{i_j}{N_j}]}(x_{1,j}) \mathbbm{1}_{[0, \frac{i_j}{N_j}]}(x_{2,j}) \int_{\frac{i_j}{N_j}}^{\frac{i_j+1}{N_j}}du\int_{\frac{i_j}{N_j}}^{\frac{i_j+1}{N_j}}du'  \partial_1K^{H'_j}(u, x_{1,j}) \\
&\hspace{6cm} \times \partial_1K^{H'_j}(u', x_{2,j})|u-u'|^{(2H'_j-2)(q-1)}\bigg)\\
&= \NN^{2\HH-1} \NN^{2-2\HH'}b_{q, \HH}^2 (\HH'(2\HH'-1))^{q-1}\\
& \times \prod_{j=1}^d\bigg(   \sum_{i_j =0}^{N_j-1} \mathbbm{1}_{[0, \frac{i_j}{N_j}]}(x_{1,j}) \mathbbm{1}_{[0, \frac{i_j}{N_j}]}(x_{2,j}) \int_{\frac{i_j}{N_j}}^{\frac{i_j+1}{N_j}}du\int_{\frac{i_j}{N_j}}^{\frac{i_j+1}{N_j}}du'  \partial_1K^{H'_j}(u, x_{1,j}) \\
&\hspace{6cm} \times \partial_1K^{H'_j}(u', x_{2,j})|u-u'|^{(2H'_j-2)(q-1)}\bigg).
\end{align*}
According to \cite[Theorem 3.2]{TudorViens1}, it is shown that for each $j$ from $1$ to $d$, the following quantity
\begin{align*}
&N_j^{2H_j-1} N_j^{2-2H'_j} \sum_{i_j=1}^{N_j-1} \mathbbm{1}_{[0, \frac{i_j}{N_j}]}(x_{1,j}) \mathbbm{1}_{[0, \frac{i_j}{N_j}]}(x_{2,j}) \int_{\frac{i_j}{N_j}}^{\frac{i_j+1}{N_j}}du\int_{\frac{i_j}{N_j}}^{\frac{i_j+1}{N_j}}du'  \partial_1K^{H'_j}(u, x_{1,j}) \\
&\hspace{6cm} \times \partial_1K^{H'_j}(u', x_{2,j})|u-u'|^{(2H'_j-2)(q-1)}
\end{align*}
converges in $L^2(\mathbb{R}^2)$  to the kernel of a standard Rosenblatt process with self-similarity $2H'_j-1$ at time $1$ (up to an explicit multiplicative constant). Since the kernel of the Rosenblatt sheet has the form of a tensor product from $1$ to $d$ of the kernel of the Rosenblatt process, (see (\ref{Def2})), 
it follows that $f_2^{\NN}$ converges to the kernel of a Rosenblatt sheet with self-similarity parameter $2\HH'-1$ evaluated at time $\textbf{1}$ up to a constant.
Therefore, the double Wiener-It\^{o} integral $\NN^{2-2\HH'}F_{2, \NN}$ in (\ref{eq:3.12}) converges in $L^2(\Omega)$ to a Rosenblatt sheet $R^{2\HH'-1}_{\textbf{1}}$ with self-similarity parameter $2\HH'-1$ evaluated at time $\textbf{1}$, which leads to the convergence of $\NN^{2-2\HH'}c_2^{-1}V_\NN$ to the same limit (up to a constant). In order to find the explicit constant, we use the fact that $\lim_{\NN \to \infty}E[(c_{1, \HH}^{-\frac12} \NN^{2-2\HH'}c_2^{-1}V_\NN)^2] = E[(R^{2\HH'-1}_{\textbf{1}})^2] =1$ to eventually obtain that $c_{1, \HH}^{-\frac12} \NN^{2-2\HH'}c_2^{-1}V_\NN$ converges in $L^2(\Omega)$ to the Rosenblatt sheet $R^{2\HH'-1}_{\textbf{1}}$ as $\NN \to \infty$ with $c_2 = q!q$.

\section*{Acknowledgements} 
I thank the authors of \cite{Reveillac3} for sharing their paper in progress with me. Another thank goes to my advisor, Prof. Ivan Nourdin, for his very careful review of the paper as well as for his comments and corrections. Finally, I deeply thank an anonymous referee for a very careful and thorough reading of this work, and for her/his
constructive remarks.

\end{document}